\documentclass{gtart_h}

\def\ifplaintex{\expandafter\ifx\csname documentclass\endcsname\relax}

\def\gtp{{\mathsurround=0pt\it $\cal G\mskip-2mu$eometry \&\ 
$\cal T\!\!$opology $\cal P\!$ublications}}  

\def\recd{{\small Received:\qua\receiveddate\ifx\reviseddate\relax
\else\qquad Revised:\qua\reviseddate\fi\par}} 

\def\lognumber#1{\def\thelognumber{#1}}
\def\volumenumber#1{\def\thevolumenumber{#1}}
\def\volumeyear#1{\def\thevolumeyear{#1}}
\def\papernumber#1{\def\thepapernumber{#1}}
\def\pagenumbers#1#2{\def\startpage{#1}\def\finishpage{#2}}
\def\published#1{\def\publishdate{#1}}

\def\received#1{\def\receiveddate{#1}}
\def\revised#1{\def\reviseddate{#1}}
\def\accepted#1{\def\accepteddate{#1}}

\long\def\asciiabstract#1{\long\def\theasciiabstract{#1}}

\let\\\par\let\thelognumber\relax\let\thevolumenumber\relax
\let\thepapernumber\relax\let\thevolumeyear\relax\let\startpage\relax
\let\finishpage\relax\let\publishdate\relax\let\receiveddate\relax
\let\reviseddate\relax\let\accepteddate\relax\let\theasciititle\relax
\let\theasciiauthors\relax
\let\theasciiabstract\relax

\let\theasciiemail\relax

\ifplaintex
\font\logobig=cmssbx10 scaled 3836
\font\logomed=cmssbx10 scaled 2557
\else
\font\logobig=cmssbx10 scaled 4200
\font\logomed=cmssbx10 scaled 2800
\fi

\long\def\makeagttitle{   
\count0=\startpage
\agt\hfill      
\hbox to 45truept{\vbox to 0pt{\vglue -13truept{\logomed A\kern -.37em{\logobig 
T}\kern -.38em G}\vss}\hss}
\break
{\small Volume \thevolumenumber\ (\thevolumeyear)
\startpage--\finishpage\nl
Published: \publishdate}

\vglue .25truein

{\parskip=0pt\leftskip 0pt plus
1fil\def\\{\par\smallskip}{\Large\bf\thetitle}\par\medskip} \vglue
0.05truein

{\parskip=0pt\leftskip 0pt plus 1fil\def\\{\par}{\sc\theauthors}
\par\medskip}%
 
\vglue 0.03truein 

{\small\leftskip 25truept\rightskip 25truept{\bf Abstract}\stdspace\theabstract

{\bf AMS Classification}\stdspace\theprimaryclass
\ifx\thesecondaryclass\relax\else; \thesecondaryclass\fi\par
{\bf Keywords}\stdspace \thekeywords\par}\vglue 7truept

}   

\ifplaintex
\font\phead=cmsl9 scaled 950
\font\pnum=cmbx10 scaled 913
\font\pfoot=cmsl9 scaled 950
\headline{\vbox to 0pt{\vskip -4.5mm\line{\small\phead\ifnum
\count0=\startpage ISSN 1472-2739 (on-line) 1472-2747 (printed)
\hfill {\pnum\folio}\else\ifodd\count0\def\\{ }%
\ifx\theshorttitle\relax\thetitle\else\theshorttitle\fi\hfill{\pnum\folio}
\else\def\\{ and }{\pnum\folio}\hfill\ifx\theshortauthors\relax\theauthors
\else\theshortauthors\fi\fi\fi}\vss}}
\footline{\vbox to 0pt{\vglue 0mm\line{\small\pfoot\ifnum\count0=\startpage
\copyright\ \gtp\hfill\else
\agt, Volume \thevolumenumber\ (\thevolumeyear)\hfill\fi}\vss}}
\else
\font=cmsl9 scaled 1050
\font\lnum=cmbx10 
\font=cmsl9 scaled 1050
\makeatletter
\def\@oddhead{{\small\lhead\ifnum\count0=\startpage ISSN 1472-2739 
(on-line) 1472-2747 (printed)\hfill {\lnum\number\count0}\else\ifodd\count0
\def\\{ }\ifx\theshorttitle\relax \thetitle \else\theshorttitle\fi\hfill
{\lnum\number\count0}\else\def\\{ and }{\lnum\number\count0}
\hfill\ifx\theshortauthors\relax 
\theauthors\else\theshortauthors\fi\fi\fi}}\def\@evenhead{\@oddhead}
\def\@oddfoot{\small\lfoot\ifnum\count0=\startpage\copyright\ \gtp\hfill\else
\agt, Volume \thevolumenumber\ (\thevolumeyear)\hfill\fi}
\def\@evenfoot{\@oddfoot}
\makeatother
\fi
\let\maketitlepage\makeagttitle
\let\makeshorttitle\maketitlepage
\let\maketitle\maketitlepage

\newwrite\gtoutfile
\long\gdef\makeheadfile{  
{\def\\{, }\def\s{ }
\immediate\openout\gtoutfile head.xxx
\immediate\write\gtoutfile{Proxy-for: \ifx\theasciiauthors\relax
\theauthors\else\theasciiauthors\fi\s<\ifx\theasciiemail\relax\theemail\else\theasciiemail\fi>}
\immediate\write\gtoutfile{\noexpand\\}
\immediate\write\gtoutfile{Authors: \ifx\theasciiauthors\relax
\theauthors\else\theasciiauthors\fi}
{\def\\{ }\immediate\write\gtoutfile{Title: \ifx\theasciititle\relax
\thetitle\else\theasciititle\fi}}
\immediate\write\gtoutfile{Subj-class: GT or SG, GR etc}
\immediate\write\gtoutfile{MSC-class: \theprimaryclass\ifx\thesecondaryclass\relax\else, \thesecondaryclass\fi}
\immediate\write\gtoutfile{Journal-ref: Algebr. Geom. Topol. \thevolumenumber\s
(\thevolumeyear) \startpage-\finishpage}
\immediate\write\gtoutfile{Comments: Published by Algebraic and
Geometric Topology at}
\immediate\write\gtoutfile{\s\s\s  http://www.maths.warwick.ac.uk/agt/AGTVol\thevolumenumber/agt-\thevolumenumber-\thepapernumber.abs.html}
\immediate\write\gtoutfile{\noexpand\\}
\immediate\write\gtoutfile{}
\ifx\theasciiabstract\relax
\immediate\write\gtoutfile{\theabstract}\else
\immediate\write\gtoutfile{\theasciiabstract}\fi
\immediate\write\gtoutfile{}
\immediate\write\gtoutfile{\noexpand\\}
\immediate\write\gtoutfile{}
\immediate\closeout\gtoutfile}}  

\def\maketitlepage{\makeagttitle\makeheadfile}
\let\makeshorttitle\maketitlepage
\let\maketitle\maketitlepage

\lognumber{46}
\volumenumber{4}
\volumeyear{2004}
\papernumber{46}
\pagenumbers{1083}{1101}
\received{4 March 2004} 
\revised{24 October 2004}
\accepted{3 November 2004}
\published{21 November 2004}

\usepackage{amsmath,amssymb,epsf}

\theoremstyle{plain}
\newtheorem{thm}{Theorem}[section]
\newtheorem{cor}[thm]{Corollary}
\newtheorem{lem}[thm]{Lemma}

\newtheorem{clm}[thm]{Claim}

\theoremstyle{definition}

\newtheorem*{rem}{Remark}
\newtheorem{qst}[thm]{Question}
\let\<\langle
\let\>\rangle

\begin{document}
\title
{Span of the Jones polynomial\\of an alternating virtual link}

\author{Naoko Kamada}

\address{Department of Mathematics,
Osaka City University,   
Sugimoto, Sumiyoshi-ku\\Osaka, 558-8585, Japan}
\email{naoko@sci.osaka-cu.ac.jp}

\begin{abstract}
For an oriented virtual link, L.H.~Kauffman defined the $f$-polynomial (Jones polynomial). 
The supporting genus of a virtual link diagram is the minimal genus of a surface in which the diagram can
 be embedded. In
this paper we show that the span of the $f$-polynomial of an alternating virtual link $L$
is determined by the number of crossings of any alternating diagram of $L$ and the 
supporting genus of the diagram. It is a generalization of Kauffman-Murasugi-Thistlethwaite's theorem. We
also prove a similar result for a virtual link diagram that is  obtained from an alternating virtual link diagram
by virtualizing one real crossing. As a consequence, such a diagram is not equivalent to a 
 classical link diagram.
\end{abstract}
\asciiabstract{%
For an oriented virtual link, L.H. Kauffman defined the f-polynomial
(Jones polynomial).  The supporting genus of a virtual link diagram is
the minimal genus of a surface in which the diagram can be
embedded. In this paper we show that the span of the f-polynomial of
an alternating virtual link L is determined by the number of crossings
of any alternating diagram of L and the supporting genus of the
diagram. It is a generalization of Kauffman-Murasugi-Thistlethwaite's
theorem. We also prove a similar result for a virtual link diagram
that is obtained from an alternating virtual link diagram by
virtualizing one real crossing. As a consequence, such a diagram is
not equivalent to a classical link diagram.}

\makeshorttitle

\section{Introduction}
\begin{figure}[htb]
\centerline{\epsfysize 6cm \epsfbox{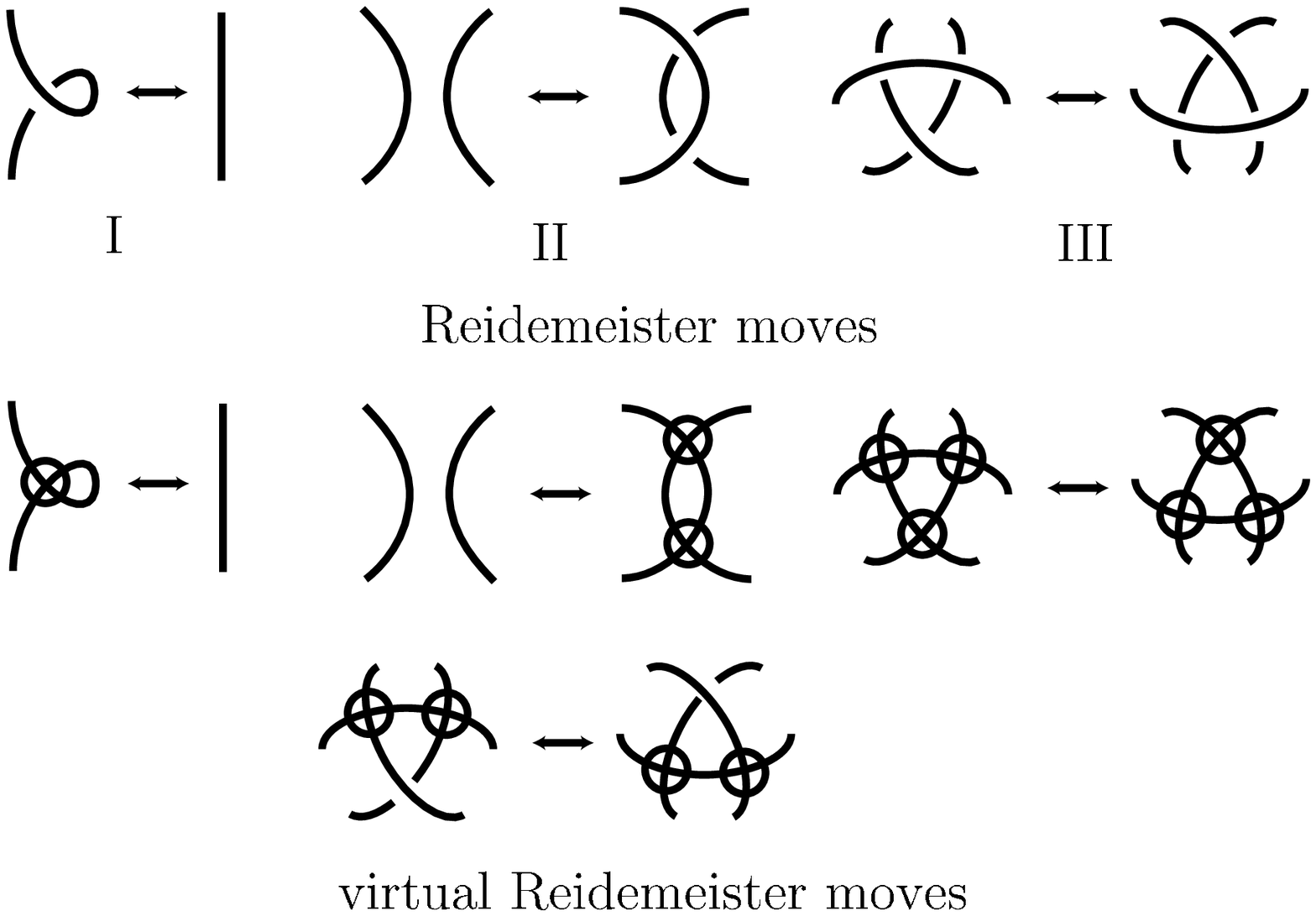}}
\caption{}\label{fig:rmove}
\end{figure}

An (oriented) {\it virtual link diagram\/} is a closed (oriented) 1-manifold  
generically immersed in ${\bf R}^2$ 
such that each double point is labeled to be either (1) a {\it real} crossing which is
indicated as usual in classical knot theory or (2) a {\it virtual} crossing which is indicated
by a small circle around the double point. 
The moves of virtual link diagrams illustrated 
in Figure~\ref{fig:rmove} are 
called {\it generalized Reidemeister moves\/}. 
Two virtual link diagrams are 
said to be {\it equivalent\/} 
if they are related by a 
finite sequence of generalized Reidemeister moves. 
A {\it virtual link\/} \cite{rGPV, rkauD} is the equivalence class of 
a virtual link diagram. Unless otherwise stated, we assume that a
virtual link is oriented.

\nocolon

Kauffman defined the $f$-polynomial $f_D(A) \in {\text{\bf Z}}[A,A^{-1}]$ of a virtual
link diagram $D$,  which is preserved under generalized Reidemeister moves, and hence it is an
invariant of a virtual link. It is also called the {\it normalized bracket polynomial\/}  or the {\it Jones
polynomial\/} 
\cite{rkauD}. 
For a virtual link $L$ represented by a virtual link diagram $D$, we define the $f$-polynomial
$f_L(A)$ of
$L$ by $f_D(A)$. The {\it span} of $f_L(A)$ is the
  maximal degree of  $f_L(A)$ minus the minimal. It is 
an invariant of a virtual link. We denote it by $\text{span}(L)$ or $\text{span}(D)$.  

By $c(D)$, we mean the number of real crossings of $D$.

\begin{thm}[Kauffman \cite{rkauA}, Murasugi \cite{rmusA}, Thistlethwaite \cite{rthi}]\label{thm:kaumura}
Let $L$ be an alternating link represented by a proper alternating connected 
link diagram
$D$. Then
we have 
$${\rm span}(L)=4c(D).$$ 
\end{thm}

Any virtual link diagram $D$  can be realized as a link diagram in a
 closed oriented surface \cite{rkauD}.
The {\it supporting genus} $g(D)$ of $D$ is the minimal genus of a closed oriented
surface in which the diagram can be realized \cite{rkk}.

Note that $g(D)$ can be calculated. Consider a link diagram ${\cal D}$ in a closed oriented surface $F$ that
realizes $D$. If some regions of the complement of ${\cal D}$ in $F$ are not open disks, replace them with
open disks. Then we obtain a link diagram realizing $D$ in a surface of genus $g(D)$. Alternatively we may
also use a formula presented in Lemma~\ref{thm:sptgenus}.

Let $D$ be a virtual link diagram. By forgetting crossing information, it is the union of immersed
circles, say $C_1,\cdots,C_{\mu}$ (for some $\mu \in {\bf N}$). The restriction of $D$ to each $C_i$ is
called a {\it component\/} of $D$, and $D$ is also called a $\mu${\it -component\/} virtual link
diagram. To state our results, we need the notion of a connected component of $D$: Consider an equivalence
relation on $C_1,\cdots , C_{\mu}$ that is the transitive closure of binary relation 
$C_i\sim C_j$ where
$C_i\sim C_j$ means that 
$C_i\cap C_j$ has at least one real crossing. Then, for an equivalence class $\{C'_1,\cdots ,C'_{\lambda}\}$, the
restriction of $D$ to $C_1'\cup\cdots\cup C_{\lambda}'$ is called a {\it connected component} of $D$. When $D$ is a connected component of itself, we say that $D$ is {\it
connected}. 

\begin{thm}\label{thm:alt2}
Let $L$ be an alternating virtual link represented by a proper
alternating virtual diagram $D$.
Then we have 
$${\rm span}(L)=4(c(D)-g(D)+m-1),$$ 
where $m$ is the number of the connected components of $D$. In particular, if $L$ is an 
alternating virtual link represented by a proper
alternating connected virtual link diagram $D$. Then
we have 
$${\rm span}(L)=4(c(D)-g(D)).$$ 
\end{thm}

Since the supporting genus of a classical link diagram is zero, Theorem~\ref{thm:alt2} 
is a generalization of Theorem~\ref{thm:kaumura}.

A similar result  was proved  in \cite{rkamC} for a link diagram in a closed oriented surface. Our argument
is essentially the same with that in \cite{rkamC}, whose basic idea is to use abstract link diagrams. 

When a virtual link diagram $D'$ is obtained from another diagram $D$ 
by replacing a real crossing $p$ of $D$ with a virtual crossing, then 
we say that $D'$ is obtained from $D$ by {\it virtualizing} the 
crossing $p$.

A virtual link diagram $D$ is said to be a {\it v-alternating} if $D$ is obtained from a 
proper alternating virtual link diagram by virtualizing one real crossing.

\begin{thm}\label{thm:valt2}
Let $D$ be a  v-alternating virtual link diagram. Then we have 
$${\rm span}(D)=4(c(D)-g(D)+m-1)+2,$$ 
where $m$ is the number of connected components of $D$. 
In particular, if $D$ is a connected
v-alternating virtual link diagram,  then $${\rm span}(D)=4(c(D)-g(D))+2.$$
\end{thm}

T. Kishino \cite{rkishi} proved that $\text{span}(D)=4c(D)-2$ when $D$ is a connected  v-alternating virtual
link diagram which is obtained from a  proper alternating classical  link diagram by virtualizing a crossing. 
His result is a special case of
Theorem~\ref{thm:valt2}, since $g(D)=1$ for such a diagram $D$ (Lemma~\ref{thm:sptgenus2}).

\begin{cor}\label{thm:valt3}
Let $D$ be a v-alternating virtual link diagram. Then
$D$ is not equivalent to a classical link diagram.
\end{cor}
\noindent
\begin{proof}
By Theorem~\ref{thm:valt2}, span$(D)$ is not a multiple of four. On the other hand, the span of the
$f$-polynomial of a classical link is a multiple of four \cite{rkauA, rmusA, rthi}. Thus we have the result. 
\end{proof}

\section{Definitions}\label{sec:def}
\begin{figure}[htb]
\centerline{\epsfysize 2.5cm \epsfbox{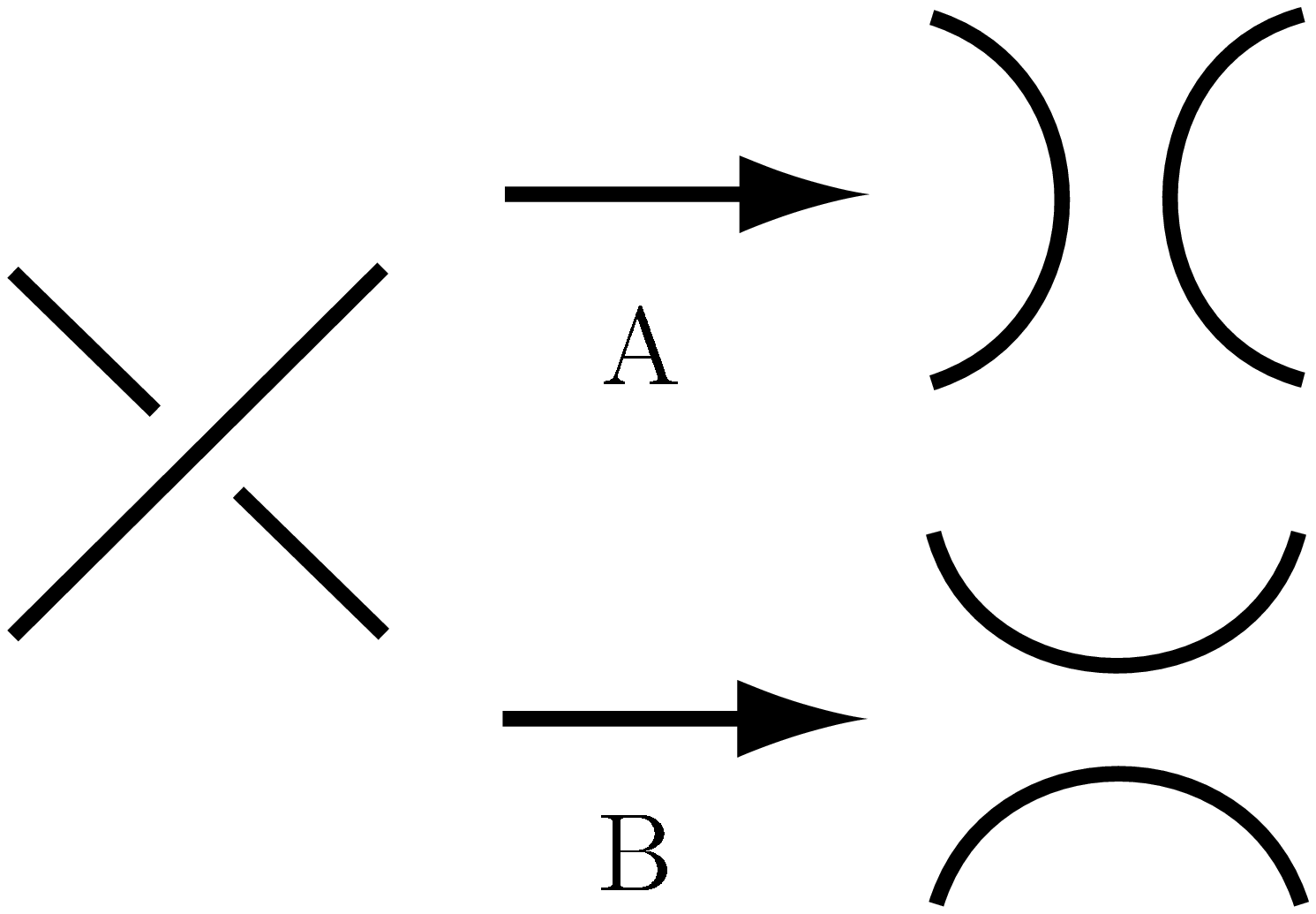}}
\caption{}\label{fig:splice}
\end{figure}

Let $D$ be an unoriented virtual link diagram. 
The replacement of the diagram in a neighborhood of a real crossing as in 
Figure~\ref{fig:splice} are called {\it A-splice\/} and {\it B-splice\/}, respectively \cite{rkauA, rkauB}.

\nocolon

A {\it state~} of $D$ 
is a virtual link diagram obtained from $D$ by doing A-splice or
B-splice at each real crossing of $D$. The {\it Kauffman bracket polynomial\/} $\<D\>$  of
$D$ is defined by
$$\<D\>=\sum_S A^{\natural (S)} (-A^2-A^{-2})^{\sharp(S)-1},$$
where $S$ runs over all states of $D$, $\natural (S)$ is the number of
A-splice minus that of B-splice used to obtain the state $S$, and
$\sharp(S)$ is the number of loops of $S$.

For an oriented virtual link diagram $D$, the {\it writhe\/} $\omega(D)$ is the number of positive
crossings minus that of negative crossings of $D$. The {\it $f$-polynomial} of $D$ is
defined by
$$f_D(A)=(-A^3)^{-\omega(D)}\<D\>.$$

\begin{thm}{\rm\cite{rkauD}}\qua
The $f$-polynomial is an invariant of a virtual link.
\end{thm}
For a virtual link $L$ represented by $D$, the $f$-polynomial $f_L(A)$ of $L$ is defined by $f_D(A)$. When
$L$ is a classical link, the $f$-polynomial $f_L(A)$ is equal to the Jones polynomial $V_L(t)$ after substituting
$A^4$ for $t$.

\begin{figure}[htb]
\centerline{\epsfysize 6cm \epsfbox{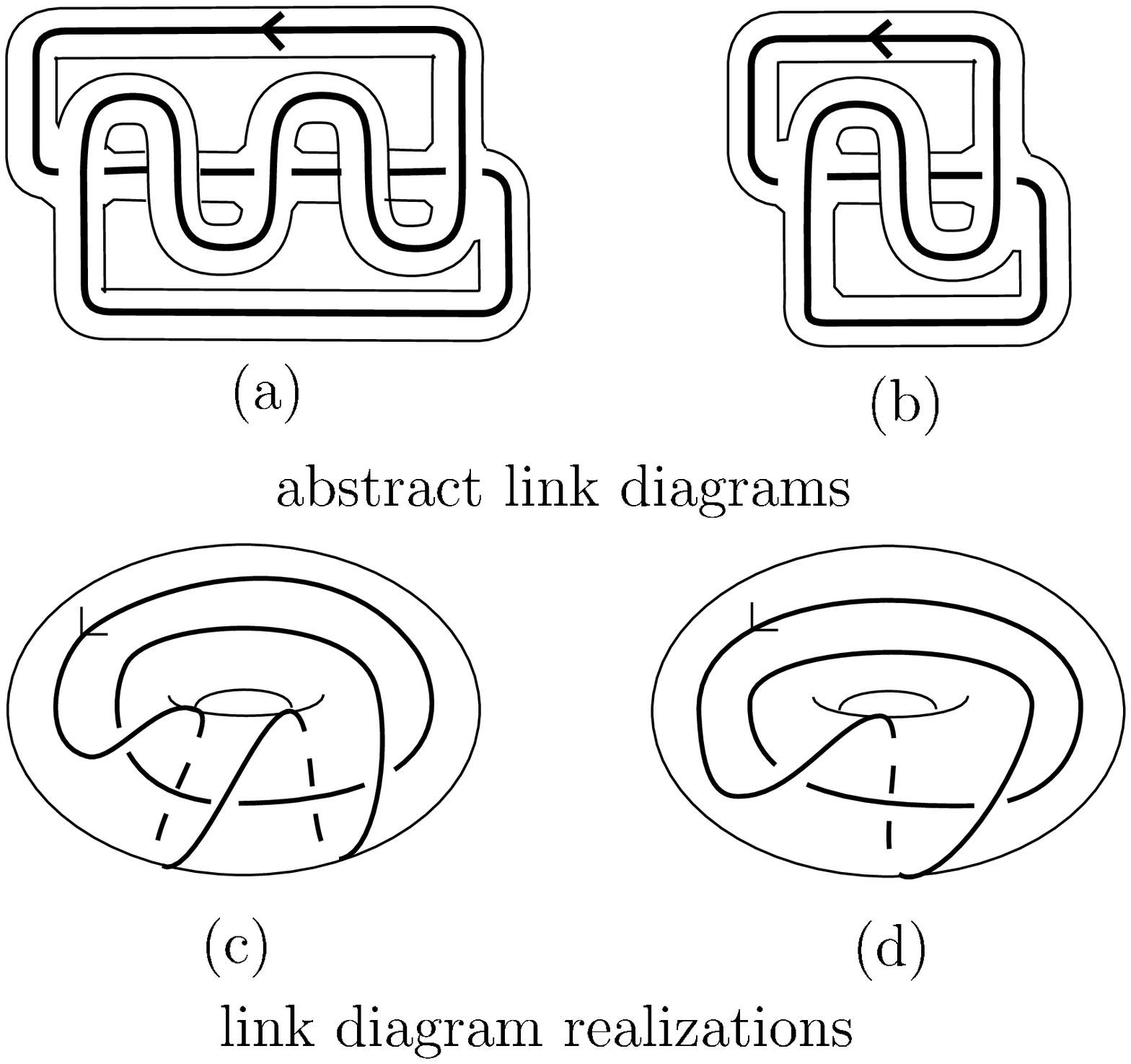}}
\caption{}\label{fig:example2}
\end{figure}

A pair $P=(\Sigma,{\cal D})$ of a compact oriented surface $\Sigma$ 
and a link diagram ${\cal D}$ in $\Sigma$ is called an {\it abstract link
diagram\/} (ALD) if $|{\cal D}|$ 
is a deformation retract of $\Sigma$, where
$|{\cal D}|$ is a graph obtained from ${\cal D}$ by replacing 
each crossing 
with a vertex. If ${\cal D}$ is an oriented link diagram, then $P$ is said to be {\it oriented\/}. Unless
otherwise stated, we assume that an ALD is oriented. If $|{\cal D}|$ is connected (or equivalently,
$\Sigma$ is connected), then $P$ is said to be {\it connected}. Two examples of  connected  ALDs are
illustrated in Figure~\ref{fig:example2} (a) and (b). 

Let  $P=(\Sigma,{\cal D})$ be an ALD.
For a closed oriented surface $F$, if there exists an embedding
$h\co \Sigma\longrightarrow F$, then $h({\cal D})$ is a link diagram in $F$. We call $h({\cal D})$ a
{\it link diagram realization\/} of  $P=(\Sigma,{\cal D})$ in $F$.
Figure~\ref{fig:example2} (c) and (d) are link diagram realizations of the ALDs illustrated in 
Figure~\ref{fig:example2} (a) and (b), respectively. 

\nocolon

The {\it supporting genus} $g(P)$ of $P=(\Sigma, {\cal D})$ is the minimal genus of a closed oriented surface
in which $\Sigma$ can be embedded \cite{rkk}.

\begin{lem}\label{thm:sptgenus}
Let  $P=(\Sigma, {\cal D})$ be an ALD, which is the disjoint union  of $m$ connected ALDs. Then
$$g(P)=\frac{2m+c({\cal D})-\sharp\partial\Sigma}{2},$$
where $c({\cal D})$ is the number of crossings of ${\cal D}$, $\partial\Sigma$ is the boundary of the
surface
$\Sigma$ and 
$\sharp\partial\Sigma$ is the number of connected components of $\partial
\Sigma$.
\end{lem}
\begin{proof}[Proof of Lemma~\ref{thm:sptgenus}]
Let $F$ be a closed oriented surface which is obtained
from $\Sigma$ by attaching $\sharp\partial\Sigma$ disks  to $\Sigma$ along the boundary $\partial\Sigma$. 
Then $g(P)=g(F)$. Since $F$ has $m$ connected components, the Euler characteristic $\chi(F)$ is
$2m-2g(F)$. On the other hand, 
$\chi(F)=\chi(\Sigma)+\sharp\partial\Sigma=\chi(|{\cal
D}|)+\sharp\partial\Sigma=-c(D)+\sharp\partial\Sigma$, since ${\cal D}$ is a 4-valent graph with
$c({\cal D})$ vertices (possibly with circle components). Thus we have the equality. 
\end{proof}

\begin{figure}[htb]
\centerline{\epsfysize 1.7cm \epsfbox{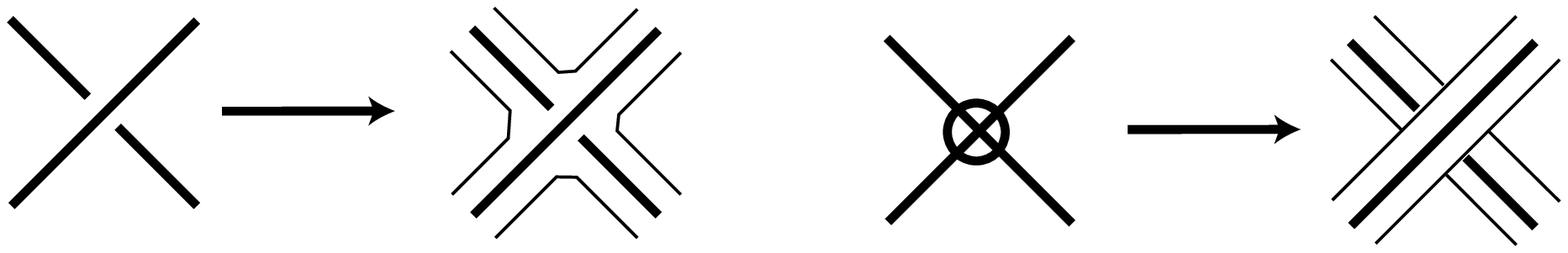}}
\caption{}\label{fig:virtualALD}
\end{figure}

\begin{figure}[htb]
\centerline{\epsfysize 2cm \epsfbox{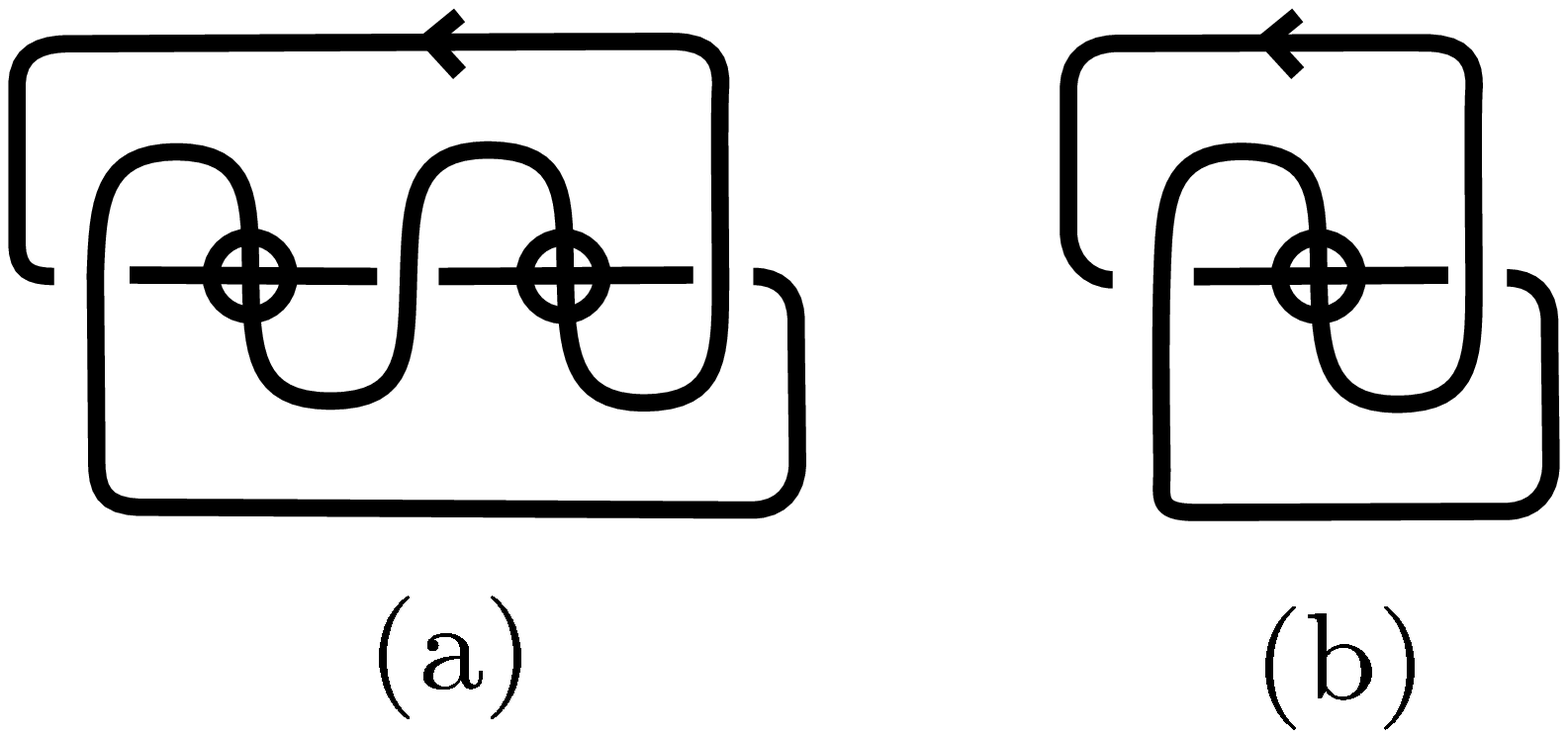}}
\caption{}\label{fig:exvir1}
\end{figure}

Let $D$ be a virtual link diagram. Consider a link diagram realization $\cal D$ of $D$ in a closed
oriented surface
$F$ and take a regular neighborhood $N({\cal D})$ of $\cal D$ in $F$. Then $(N({\cal D}),{\cal D})$ is an ALD. 
We call it the {\it ALD associated with} $D$, and denote it by $\phi(D)$. 
(Note that $\phi(D)$, up to homeomorphism, does not depend on the choice of $F$ and the realization $\cal D$
in $F$.)  An easy method to obtain
$\phi(D)$ is illustrated in  Figure~\ref{fig:virtualALD} (see
\cite{rkk} for details).  For example, the ALDs  illustrated in Figure~\ref{fig:example2}
(a) and (b) are the ALDs associated with the virtual link diagrams in Figure~\ref{fig:exvir1}
(a) and (b), respectively. 

\nocolon

\nocolon

\begin{lem}
Let $D$ be a virtual link diagram and let $\phi(D)=P=(\Sigma, {\cal D})$ be the ALD associated with
$D$. Then we have
\begin{itemize}
\item[{\rm (1)}]
$g(P)=g(D)$
\item[{\rm (2)}]
$P$ is connected if and only if $D$ is connected.
\end{itemize}
\end{lem}

\begin{proof}
It is obvious from the definition. 
\end{proof}

\begin{rem}
Let $P=(\Sigma, {\cal D})$  and $P'=(\Sigma', {\cal D}')$ be ALDs. We say that $P'$ is obtained from $P$ by an
{\it abstract Reidemeister move\/} if  there are embeddings 
$h\co \Sigma\longrightarrow F$ and
$h'\co \Sigma'\longrightarrow F$  into a closed oriented surface $F$ such that the  link diagram $h({\cal D}')$ 
is obtained from
$h({\cal D})$ by a Reidemeister move in $F$. Two ALDs $P=(\Sigma, {\cal D})$  and $P'=(\Sigma', {\cal
D}')$ are {\it equivalent} if there exists a finite sequence of ALDs, $P_0,
P_1,\cdots ,P_u$,  with $P_0=P$ and 
$P_u=P'$ such that
$P_{i+1}$ is obtained from $P_i$ by an abstract Reidemeister move. 
An {\it abstract link\/} is such an equivalence class (cf. \cite{rkk}). It is proved in
\cite{rkk} that two
virtual link diagrams $D$ and $D'$ are equivalent if and only if the associated ALDs, $\phi(D)$ and $\phi(D')$,
 are equivalent; namely, the map 
$$\phi\co \{\text{\rm virtual link diagrams}\}\longrightarrow \{\text{abstract link diagrams}\}$$
induces a bijection
$$\{\text{\rm virtual links}\}\longrightarrow \{\text{abstract links}\}.$$
\end{rem}

\begin{figure}[htb]
\centerline{\epsfysize 4cm \epsfbox{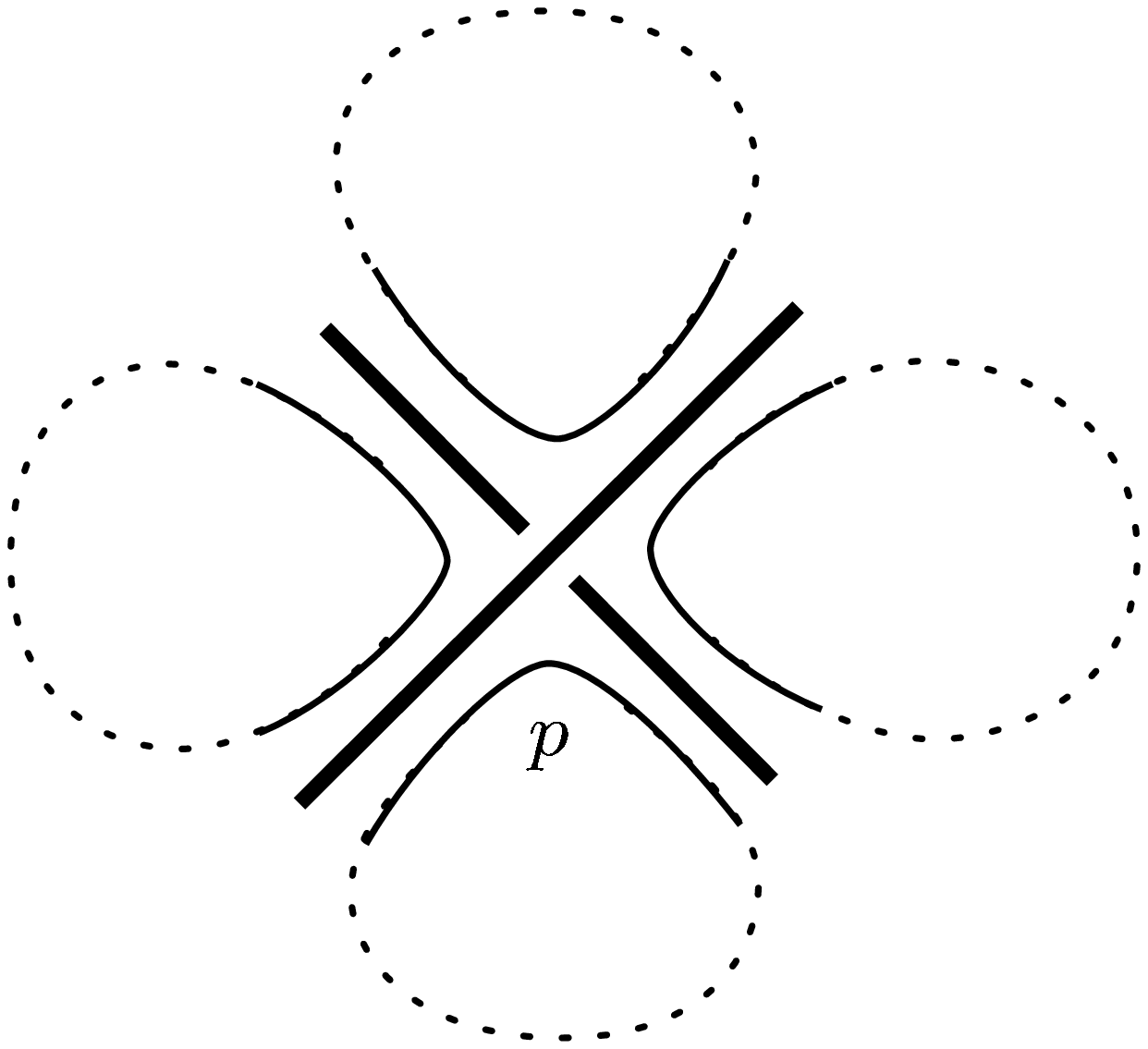}}
\caption{}\label{fig:prpvkd}
\end{figure}

Let $P=(\Sigma, {\cal D})$ be an ALD. A crossing $p$ of ${\cal D}$ is {\it proper} if four connected
components of 
$\partial\Sigma$  passing through the neighborhood of $p$ are all distinct. See Figure~\ref{fig:prpvkd}.
When every crossing of $\cal D$ is proper, we say that $P$ is {\it proper}.
Let $D$ be a virtual link diagram and $\phi(D)=(\Sigma,{\cal D})$ the ALD associated
with $D$. A real crossing of $D$ is said to be {\it proper} if  the  corresponding crossing of
${\cal D}$ is proper. 
A virtual link diagram $D$ is said to be {\it proper\/} if each crossing of $D$ is proper (or equivalently if
$\phi(D)$ is a proper ALD).

\nocolon

\begin{figure}[htb]
\centerline{\epsfysize 4cm \epsfbox{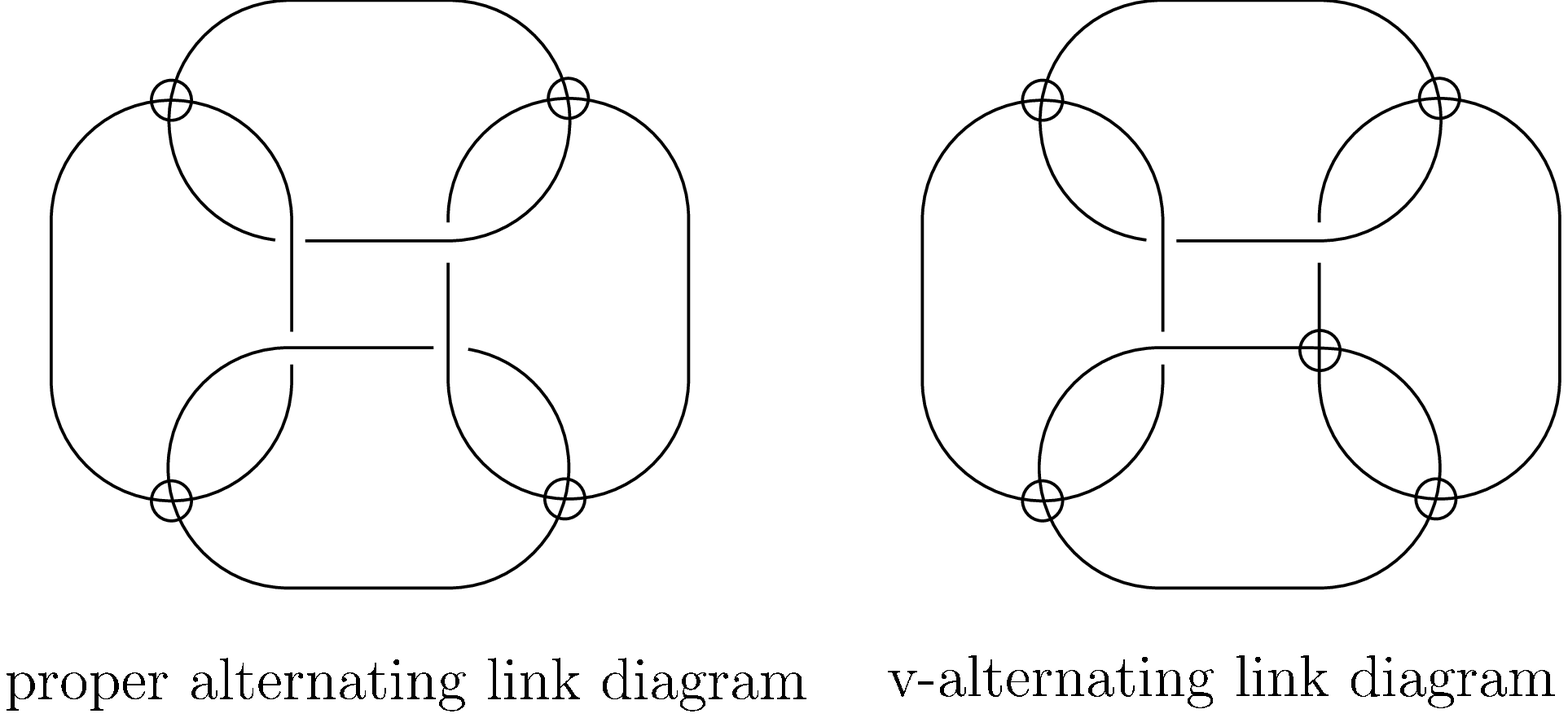}}
\caption{}\label{fig:expaltvalt1}
\end{figure}

The left hand side of  Figure~\ref{fig:expaltvalt1} is a proper alternating virtual link diagram and the right hand
side is a non-proper virtual link diagram. The right hand side is a 
v-alternating virtual link diagram obtained from the left diagram by virtualizing a real crossing.

\nocolon

\section{Checkerboard coloring}

Let $P=(\Sigma, {\cal D})$ be an ALD. We say that $P$ is {\it chekerboard colorable} if we can assign two
colors (black and white) to the region of $\Sigma\setminus|{\cal D}|$ such that two adjacent regions with an
arc of $|{\cal D}|$ have distinct colors, where $|{\cal D}|$ is the graph obtained from $\cal D$ by
assuming each crossing to be a vertex of degree four. A {\it checkerboard coloring} of $P$ is such
an assignment of colors.  

\begin{figure}[htb]
\centerline{\epsfysize 2.5cm \epsfbox{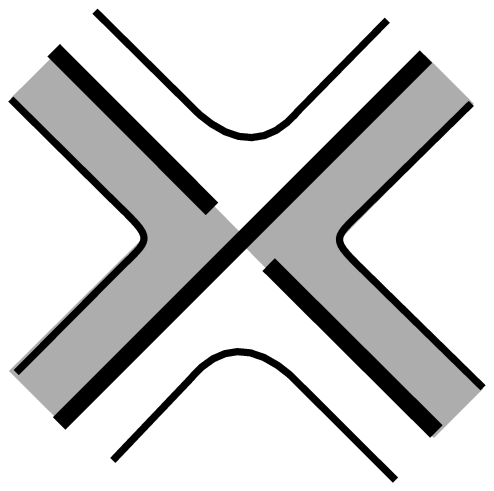}}
\caption{}\label{fig:chkbd1}
\end{figure}

If $P$ is an alternating ALD, then it has a checkerboard coloring such that for each crossing, the
regions around each crossing are colored as in Figure~\ref{fig:chkbd1}. (This fact is seen as follows:
 Walk on any knot component of $\cal D$ and look at the right hand side. When we pass a crossing as an
over-arc, or as an under-arc, the right is colored black, or white respectively. Since $\cal D$ is alternating, we
have a coherent coloring.) We call such a coloring a {\it canonical checkerboard coloring} of an alternating 
ALD, which is unique unless
$P$ has a connected component without crossings. 

\nocolon

\begin{figure}[htb]
\centerline{\epsfysize 2.5cm \epsfbox{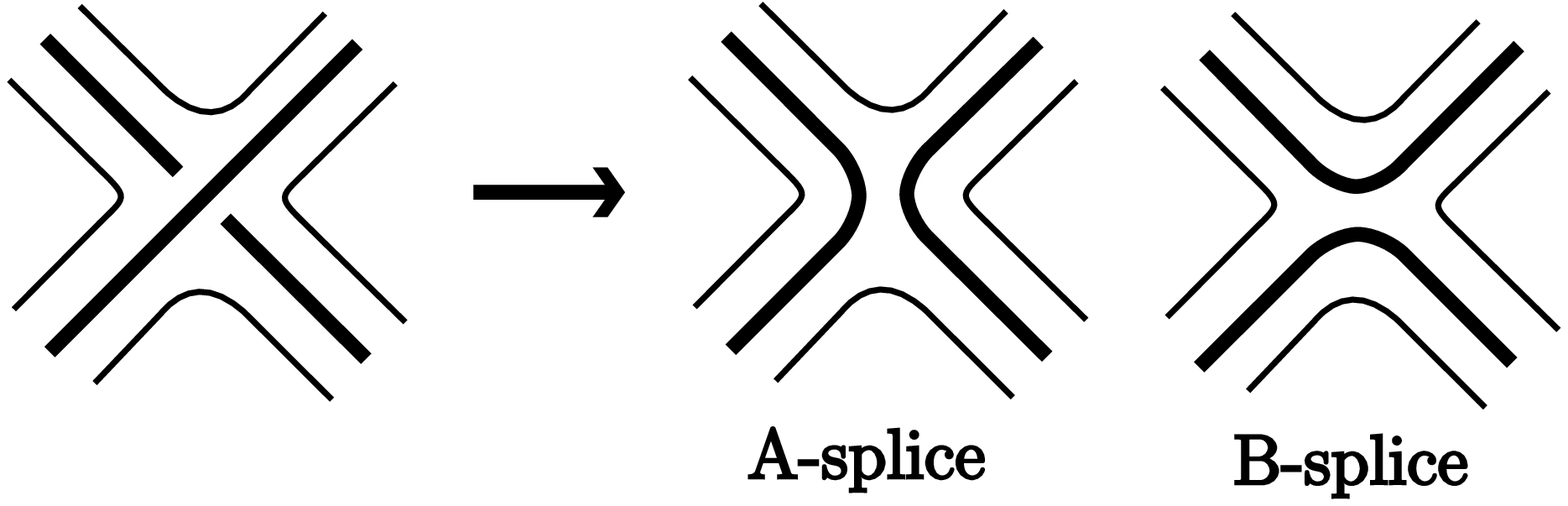}}
\caption{}\label{fig:chkbd2}
\end{figure}

Let $P=(\Sigma, {\cal D})$ be an ALD and let ${\cal S}_A$ (or ${\cal S}_B$, resp.) be the state of $\cal D$
obtained from $\cal D$ by doing A-splice (resp. B-splice) for every crossing. 
(See Figure~\ref{fig:chkbd2}. The states on $\Sigma$ are no longer ALDs.)

\nocolon

\begin{figure}[htb]
\centerline{\epsfysize 2.8cm \epsfbox{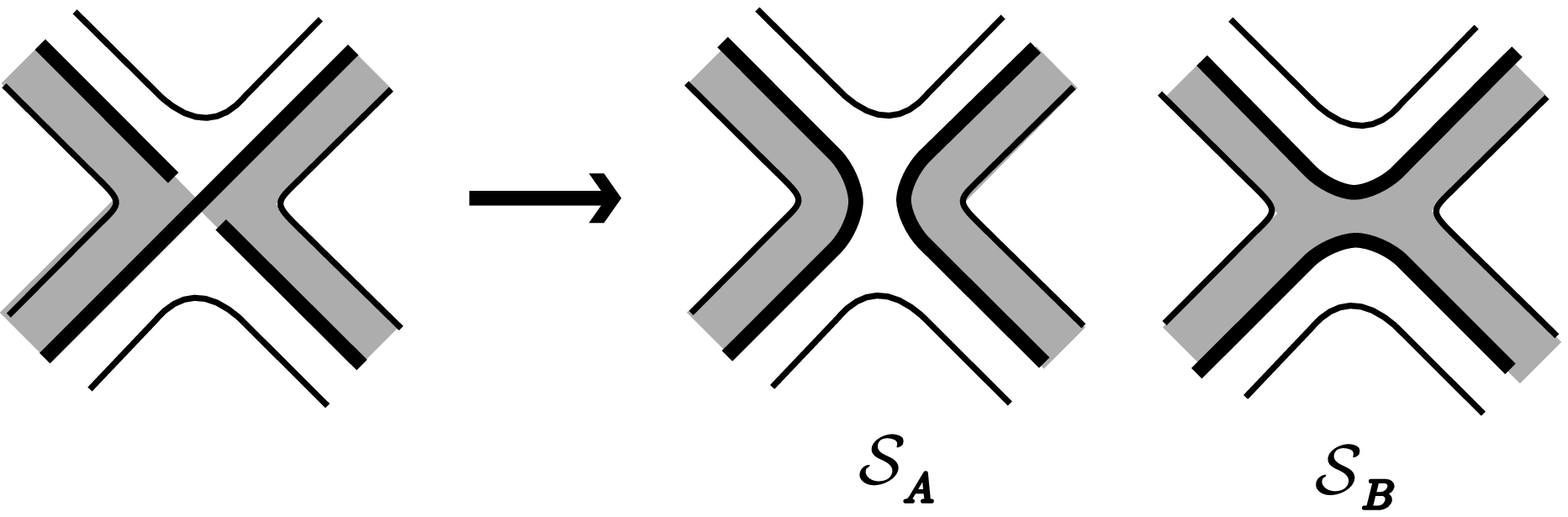}}
\caption{}\label{fig:chkbd3}
\end{figure}

Suppose that  $P=(\Sigma, {\cal D})$  be alternating, and consider a canonical checkerboard coloring of $P$.
Then $(\Sigma, {\cal S}_A)$ and $(\Sigma, {\cal S}_B)$ inherit checkerboard colorings. See
Figure~\ref{fig:chkbd3}. Black regions of $(\Sigma, {\cal S}_A)$ are annuli. 
Thus we have a
one-to-one correspondence
$$\{\text{\rm the loops of } {\cal S}_A\}\longrightarrow \{\text{\rm the loops of }\partial\Sigma
\text{
\rm in black regions}\}$$
so that a loop of ${\cal S}_A$ and the corresponding loop of $\partial \Sigma$ bound an annulus colored black. 
Similarly, white regions of $(\Sigma, {\cal S}_B)$ are annuli. Thus we have a one-to-one correspondence 
$$\{\text{\rm the loops of } {\cal S}_B\}\longrightarrow \{\text{\rm the loops of }\partial\Sigma\text{
\rm in white regions}\}$$ 
so that a loop of ${\cal S}_B$ and the corresponding loop of $\partial \Sigma$
bound an annulus colored white. Thus we have the following.

\nocolon

\begin{figure}[htb]
\centerline{\epsfysize 9.5cm \epsfbox{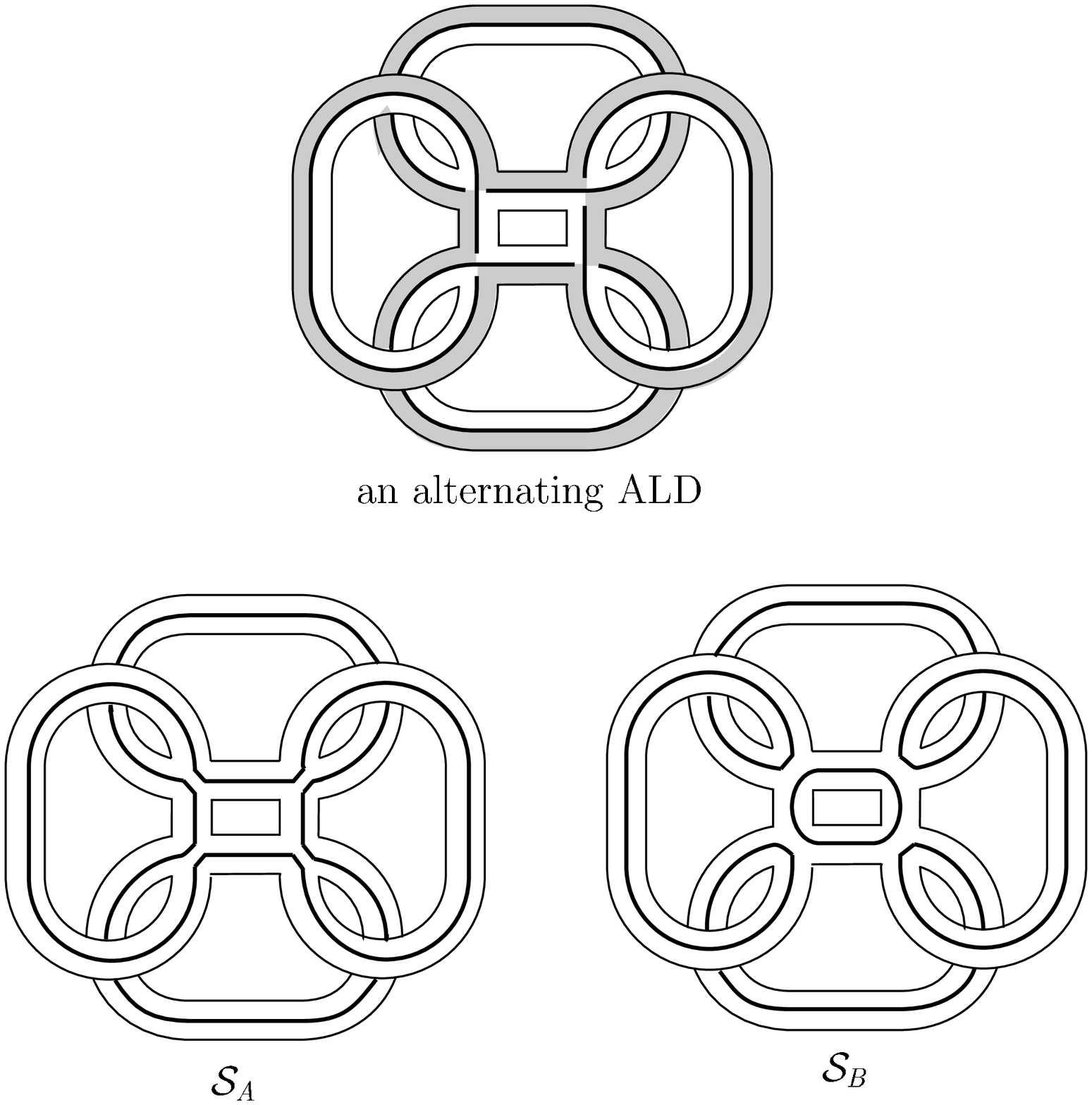}}
\caption{}\label{fig:altstate1}
\end{figure}

\begin{lem}\label{thm:stateloop}
In the situation above, there is a bijection
$$
\{\text{\rm the loops of }{\cal S}_A\}\cup \{\text{\rm the loops of
}{\cal S}_B\}\longrightarrow\{\text{\rm the loops of }\partial\Sigma\}
.$$
\end{lem}
We have an example of an
alternating ALD with a canonical checkerboard coloring and the states ${\cal S}_A$ and ${\cal S}_B$ in
Figure~\ref{fig:altstate1}.   

\nocolon

\begin{lem}\label{thm:state1}
Let $P=(\Sigma, {\cal D})$ be an alternating ALD, and let  ${\cal S}_A$ (or ${\cal S}_B$, resp.) be the state of
${\cal D}$ obtained from $\cal D$ by doing A-splice (resp. B-splice) for every crossing. For a crossing $p$ of
$\cal D$, let
$l_1(p)$ and
$l_2(p)$ be the loops of
${\cal S}_A$ (or
$l'_1(p)$ and $l'_2(p)$ be the loops of ${\cal S_B}$) that pass through the neighborhood of $p$. If $p$ is a
proper crossing, then $l_1(p)\ne l_2(p)$ and $l'_1(p)\ne l'_2(p)$. 
\end{lem}
\begin{proof}
Since $p$ is a proper crossing, the four loops of $\partial\Sigma$ appearing around $p$ arc all distinct.
Since $P$ is alternating, it has a canonical checkerboard coloring and there is a one-to-one
correspondence as in Lemma~\ref{thm:stateloop}. Then $l_1(p)$,
$l_2(p)$, $l'_1(p)$ and $l'_2(p)$ correspond to the four distinct loops of $\partial\Sigma$ around
$p$. Thus $l_1(p)\ne l_2(p)$ and  $l'_1(p)\ne l'_2(p)$. 
\end{proof}

\section{Proofs of Theorems~\ref{thm:alt2} and \ref{thm:valt2}}\label{sec:prf}

We denote  the maximal  (or  minimal, resp.) degree
of a Laurent polynomial $\eta$ by ${\rm maxd}(\eta)$ (resp.
${\rm mind}(\eta)$).  For a state $S$ of a virtual link diagram $D$, let $\<S|D\>$ stand for $A^{\natural
S}(-A^2-A^{-2})^{\sharp S-1}$. 

\begin{proof}[Proof of Theorem~\ref{thm:alt2}]
Let $D$ be a proper alternating virtual link diagram of $m$ connected components, and let
$P=(\Sigma,{\cal D})$ be the ALD associated with $D$.  
Let $S_A$ (or $S_B$ resp.) be the state of $D$ obtained from $D$ by doing A-splice (resp.
B-splice) at each crossing of $D$, and let ${\cal S}_A$ (resp. ${\cal S}_B$) be the corresponding state of $\cal
D$ in
$\Sigma$. 

Let $S_A(j)$ (or $S_B(j)$,  resp.) be a state obtained from $S_A$ (resp. $S_B$) by changing
A-splices (resp. B-splices) to B-splices (resp. A-splices) at $j$ crossings of $D$. 

\begin{clm}\label{thm:statenum1}
$\sharp S_A(1)=\sharp S_A-1$ and $\sharp S_B(1)=\sharp S_B-1$.
\end{clm}
\begin{proof}
Let $S_A(1)$ be obtained from $S_A$ by changing A-splice to B-splice at a crossing point $\tilde{p}$ of $D$.
Let ${\cal S}_A(1)$  be the corresponding state of ${\cal D}$, and let $p$ be the crossing
of ${\cal D}$ corresponding to $\tilde{p}$. Since $D$ is proper, the crossing $p$ is proper. We prove the
former equality for the corresponding ALD version; namely, $\sharp {\cal S}_A(1)=\sharp {\cal S}_A-1$. In
the situation of Lemma~\ref{thm:state1},
$l_1(p)\ne l_2(p)$. Since ${\cal S}_A(1)$ is obtained from ${\cal S}_A$ by changing A-splice with
B-splice at
$p$, two distinct loops $l_1(p)$ and $l_2(p)$ become a single loop. Hence $\sharp{\cal S}_A(1)={\cal
S}_A-1$. Therefore we have $\sharp S_A(1)=\sharp S_A-1$. Similarly, we have $\sharp S_B(1)=\sharp
S_B-1$. 
\end{proof}

\begin{clm}\label{thm:statenum2}
$\sharp S_A(j)\leq\sharp S_A+j-2$ and $\sharp S_B(j)\leq\sharp S_B+j-2$ for $j=1,\cdots,c(D)$.
\end{clm}
\begin{proof}
Any $S_A(k)$ , $k=1,\cdots,c(D)$, is obtained from some $S_A(k-1)$ by changing A-splice to B-splice at a
crossing. Then 
$$\sharp S_A(k-1)-1\leq \sharp S_A(k)\leq \sharp S_A(k-1)+1.$$
Thus $\sharp S_A(j)\leq \sharp S_A(1)+j-1$. By Claim~\ref{thm:statenum1}, we have $\sharp S_A(j)\leq
 \sharp S_A+j-2$. Similarly, we have $\sharp S_B(j)\leq \sharp S_B+j-2$. 
\end{proof}

Now we continue the proof of Theorem~\ref{thm:alt2}. By definition,  
\begin{eqnarray}
{\rm maxd}(\<S_A|D\>)&=&{\rm maxd}(A^{c(D)}(-A^2-A^{-2})^{\sharp S_A-1})\nonumber\\
&=&c(D)+2\sharp S_A-2 \label{eqn:astate1}
\end{eqnarray}
and 
\begin{eqnarray}
{\rm mind}(\<S_B|D\>)&=&{\rm mind}(A^{-c(D)}(-A^2-A^{-2})^{\sharp S_B-1})\nonumber\\
&=&-c(D)-2\sharp S_B+2.\label{eqn:bstate1}
\end{eqnarray}
For a state $S_A(j)$ for $j=1,\cdots,c(D)$,  using Claim~\ref{thm:statenum2}, we have
\begin{eqnarray}
{\rm maxd}(\<S_A(j)|D\>)&=&{\rm maxd}(A^{{c(D)-2j}}(-A^2-A^{-2})^{\sharp
S_A(j)-1})\nonumber\\
&=&c(D)-2j+2\sharp S_A(j)-2. \nonumber\\
&\leq& c(D)+2\sharp S_A-6.
\label{eqn:astate2}
\end{eqnarray}
For a state $S_B(j)$ for $j=1,\cdots,c(D)$, using Claim~\ref{thm:statenum2}, we have
\begin{eqnarray}
{\rm mind}(\<S_B(j)|D\>)&=&{\rm mind}(A^{-c(D)+2j}(-A^2-A^{-2})^{\sharp
S_B(j)-1})\nonumber\\
&=&-c(D)+2j-2\sharp S_B(j)+2.\nonumber\\
&\geq& -c(D)-2\sharp S_B+6.
\label{eqn:bstate2}
\end{eqnarray}
From (\ref{eqn:astate1}), (\ref{eqn:bstate1}), (\ref{eqn:astate2}), (\ref{eqn:bstate2}) we have
$$\left\{\begin{array}{l}
{\rm maxd}(\<D\>)=c(D)+2\sharp S_A-2,\\
{\rm mind}(\<D\>)=-c(D)-2\sharp S_B+2.
\end{array}\right.$$
Thus 
$${\rm span}(D)=2c(D)+2(\sharp S_A+\sharp S_B)-4.$$
By Lemma~\ref{thm:stateloop}, we have $\sharp S_A+\sharp S_B=\sharp
{\cal S}_A+\sharp {\cal S}_B=\sharp\partial\Sigma$. Therefore
$${\rm span}(D)=2c(D)+2\sharp \partial \Sigma-4.$$
By Lemma~\ref{thm:sptgenus}, we have the desired equality. 
\end{proof}

\begin{figure}[htb]
\centerline{\epsfysize 3cm \epsfbox{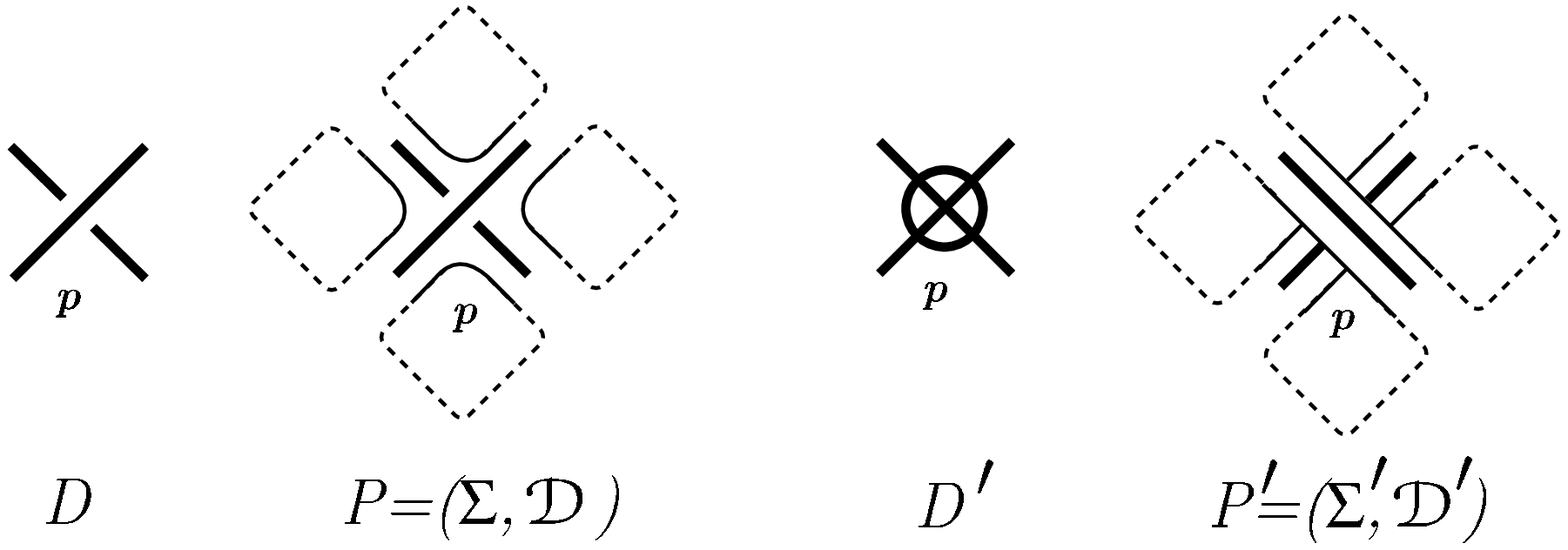}}
\caption{}\label{fig:valtstate1}
\end{figure}

\begin{proof}[Proof of Theorem~\ref{thm:valt2}]
Let $D'$ be a v-alternating virtual link diagram obtained from a proper alternating virtual link diagram $D$
 by virtualizing a real crossing
$p$ of
$D$, and let $P'=(\Sigma', {\cal D}')$ be the ALD associated with $D'$. Note that $P'=(\Sigma', {\cal D}')$
is obtained from the ALD, $P=(\Sigma, {\cal D})$, associated with $D$ by changing the
neighborhood of the crossing which corresponds to $p$ of $D$ as in Figure~\ref{fig:altvalt1}.

\nocolon
\begin{figure}[htb]
\centerline{\epsfysize 2cm \epsfbox{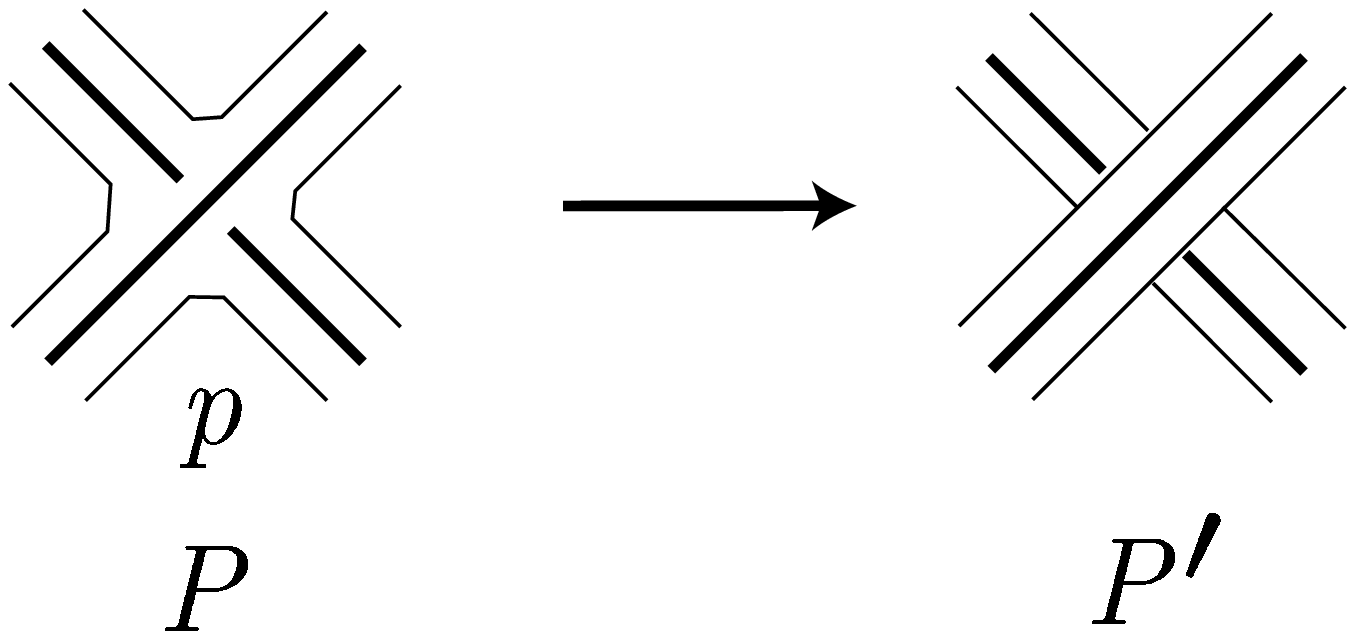}}
\caption{}\label{fig:altvalt1}
\end{figure}

Let $S_A$ (or $S_B$ resp.) be the state of $D$ obtained by doing A-splice (resp. B-splice) at each crossing, and
let $S'_A$ (resp. $S'_B$) be the state of $D'$ obtained by doing A-splice (resp.
B-splice) at each crossing. 
$S'_A$ (or $S'_B$ resp.) is obtained from $S_A$ (resp. $S_B$ ) by connecting two connected components of 
$S_A$ which pass through the neighborhood of $p$ as in Figure~\ref{fig:valtstate3}.

\nocolon
\begin{figure}[htb]
\centerline{\epsfysize 3.7cm \epsfbox{ 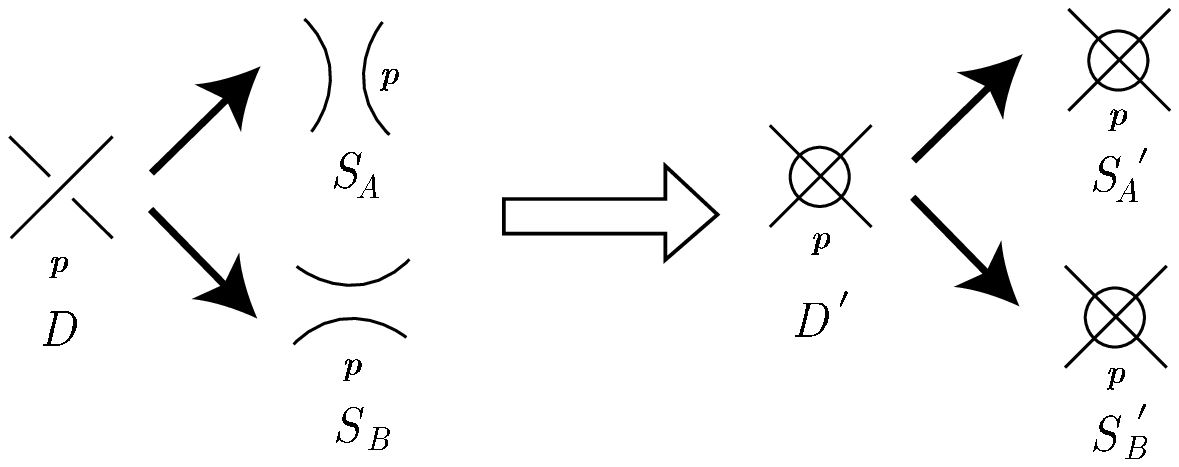}}
\caption{}\label{fig:valtstate3}
\end{figure}

Let $S'_A(j)$ (or $S'_B(j)$, resp.) be a state obtained from $S'_A$ (resp. $S'_B$) by changing
A-splices (resp. B-splices) to B-splices (resp. A-splices) at $j$ crossings of $D'$. 

\begin{clm}\label{thm:statenum3}
\begin{itemize}
\item[{\rm (1)}]
$\sharp S'_A-1\leq \sharp S'_A(1) \leq \sharp S'_A$ and 
$\sharp S'_B-1\leq \sharp S'_B(1) \leq \sharp S'_B$.
\item[{\rm (2)}]
$\sharp S'_A(j)\leq \sharp S'_A+j-1$ and 
$\sharp S'_B(j)\leq \sharp S'_B+j-1$ for $j=1,2,\cdots,c(D')$.
\end{itemize}
\end{clm}
\begin{proof}
Any $S'_A(k)$, $k=1,\cdots,c(D')$, is obtained from some $S'_A(k-1)$ by changing A-splice to B-splice at
a crossing. Then 
\begin{equation}
\sharp S'_A(k-1)-1\leq \sharp S'_A(k)\leq \sharp S'_A(k-1)+1.
\label{eqn:statevalt1}
\end{equation}
In particular, $\sharp S'_A-1\leq \sharp S'_A(1)\leq \sharp S'_A+1$. If $\sharp S'_A(1)=\sharp S'_A+1$, then
$\sharp S_A(1)=\sharp S_A+1$ (see  Figure~\ref{fig:valtstate2}). It contradicts that  $D$ is proper (recall
Claim~\ref{thm:statenum1}). Thus we have $\sharp S'_A-1 \leq \sharp S'_A(1)\leq \sharp S'_A$.  By
(\ref{eqn:statevalt1}), $\sharp S'_A(j)\leq \sharp S'_A(1)+j-1$. Hence $\sharp S'_A(j)\leq\sharp S'_A+j-1$.
Similarly we have $\sharp S'_B-1\leq \sharp S'_B(1)\sharp S'_B$ and $\sharp S'_B(j)\leq \sharp
S'_B+j-1$.
\end{proof}

\nocolon
\begin{figure}[htb]
\centerline{\epsfysize 8cm \epsfbox{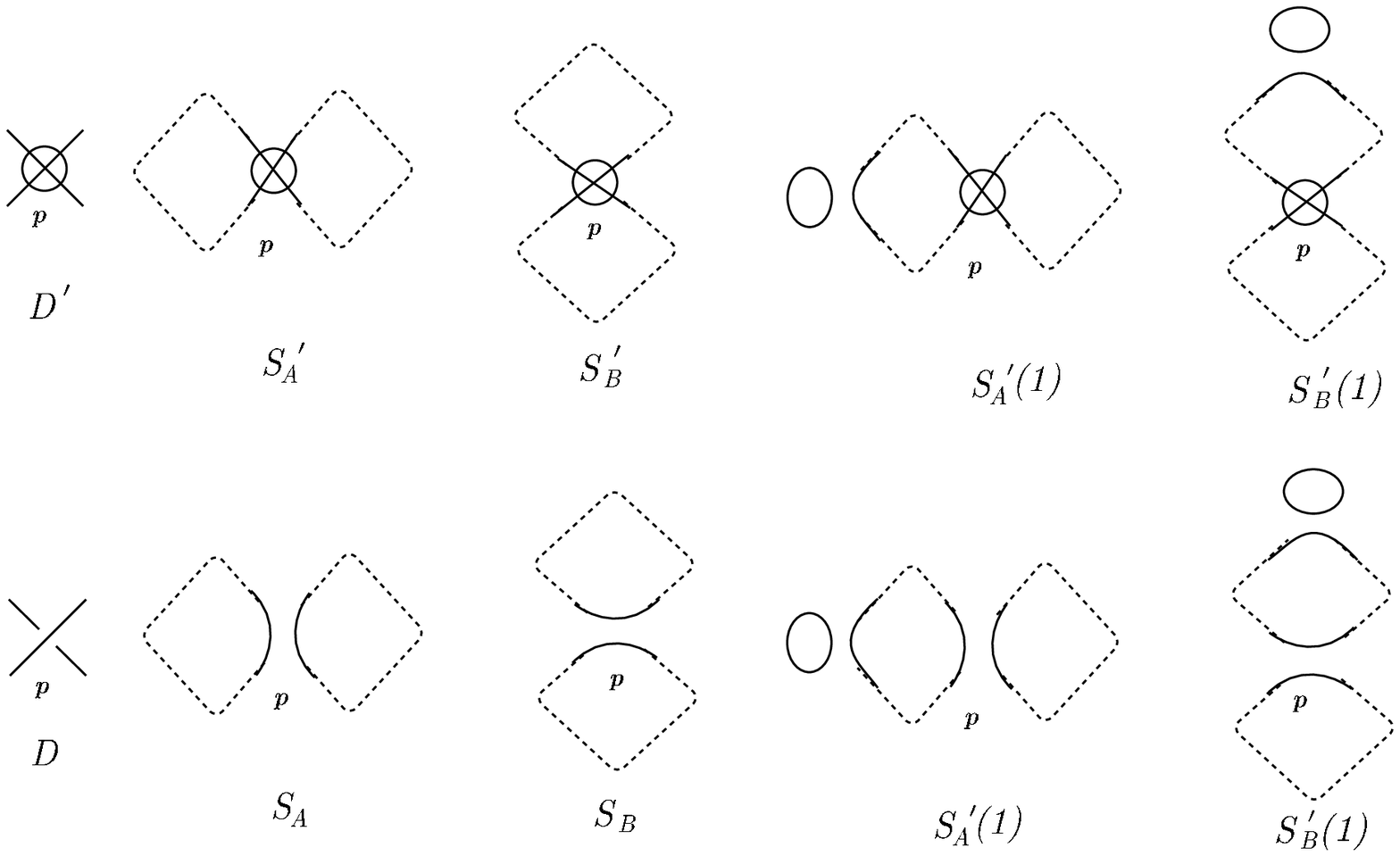}}
\caption{}\label{fig:valtstate2}
\end{figure}

By definition, we have 
$${\rm maxd}(\<S'_A|D'\>)=c(D')+2\sharp S'_A-2$$
and 
$${\rm mind}(\<S'_B|D'\>)=-c(D')-2\sharp S'_B+2.$$
For a state $S'_A(j)$ and $S'_B(j)$, using Claim~\ref{thm:statenum3}, we have
$$\begin{array}{ll}
{\rm maxd}(\<S'_A(j)|D'\>)&=c(D')-2j+2\sharp S'_A(j)-2\\
&\leq c(D')+2\sharp S'_A-4
\end{array}$$
and
$$\begin{array}{ll}
{\rm mind}(\<S'_B(j)|D'\>)&=-c(D')+2j-2\sharp S'_B(j)+2
\\
&\geq -c(D')-2\sharp S'_B+4.
\end{array}
$$  
Therefore, we have 
$$\left\{
\begin{array}{l}
{\rm maxd}\<D'\>=c(D')+2\sharp S'_A-2\\
{\rm minxd}\<D'\>=-c(D')-2\sharp S'_B+2
\end{array}\right.
$$
and
$${\rm span}(D')=2c(D')+2(\sharp S'_A+\sharp S'_B)-4.$$
Since $p$ is proper, by Lemma~\ref{thm:state1}, we see that $\sharp S'_A=\sharp S_A-1$ and $\sharp
S'_B=\sharp S_B-1$. By Lemma~\ref{thm:stateloop}, we have ${\rm span}(D')=2c(D')+2(\sharp S_A+\sharp
S_B)-8=2c(D')+2\sharp \partial\Sigma-8$. 

\begin{clm}\label{thm:bdynum1}
$\sharp\partial\Sigma'=\sharp\partial\Sigma-3$.
\end{clm}
\begin{proof}
Since $p$ is a proper crossing, the four loops of $\partial\Sigma$ around $p$ are all distinct.
After changing $P=(\Sigma, {\cal D})$ to $P'=(\Sigma',{\cal D}')$ as in Figure~\ref{fig:altvalt1}, the four
loops become a single loop of $\partial\Sigma'$
(see Figure~\ref{fig:valtstate1}). 
\end{proof}

Thus ${\rm span}(D')=2c(D')+2\sharp\partial\Sigma'-2$.
By Lemma~\ref{thm:sptgenus}, we have
$g(D')=(2m+c(D')-\sharp \partial\Sigma')/2$.
Therefore ${\rm span}(D')=4(c(D')-g(D')+m-1)+2$. This completes the proof of Theorem~\ref{thm:valt2}.
\end{proof}

\nocolon

\begin{lem}\label{thm:sptgenus2}
Suppose that a virtual link diagram $D'$  is obtained from
a virtual link diagram  $D$ by virtualizing a crossing $p$ of $D$.
If $p$ is proper, then $g(D') = g(D) +1$.
\end{lem}

\begin{proof}
Let $P =(\Sigma, {\cal D})$ and $P' =(\Sigma', {\cal D}')$ be the
ALDs associated with $D$ and $D'$.
Since $p$ is proper, the numbers of connected components of
$\Sigma$ and $\Sigma'$ must be the same, and as we saw in Claim~\ref{thm:bdynum1}
(Figure~\ref{fig:valtstate1}),
$\sharp \partial \Sigma'= \sharp \partial \Sigma -3$.
Since $c(D') = c(D) -1$, by Lemma 8, we seen that $g(\Sigma') = g(\Sigma)+1$.
Thus $g(D') = g(D) +1$. 
\end{proof}

\section{2-braid virtual link}
\begin{figure}[htb]
\centerline{\epsfysize 3.5cm \epsfbox{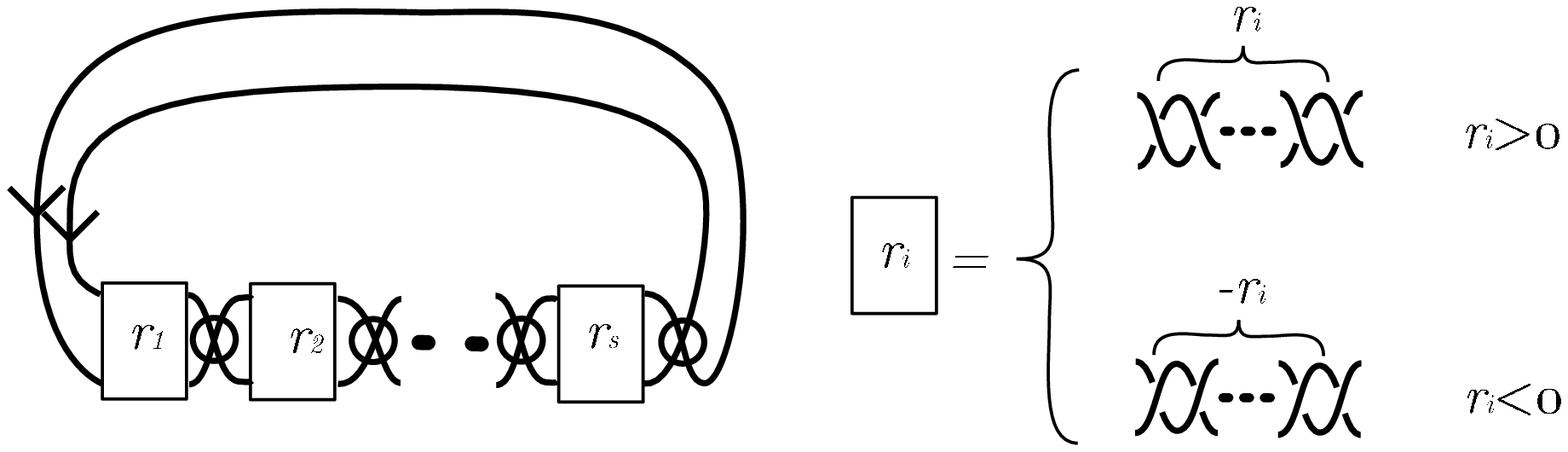}}
\caption{}\label{fig:2brknot}
\end{figure}

For non-zero integer $r_1,\cdots,r_s$, we denote by $K(r_1,\cdots , r_s)$ a  virtual link diagram illustrated  in
Figure~\ref{fig:2brknot}. The virtual link represented by this diagram is also denoted by $K(r_1,\cdots,r_s)$.
M. Murai \cite{rmurai} proved that $K(r_1)$ and $K(r_1, r_2) $ are not  classical links and that 
$K(r_1)$ and $K(r_2, r_3)$ are distinct virtual links.

\nocolon

\begin{figure}[htb]
\centerline{\epsfysize 1.3cm \epsfbox{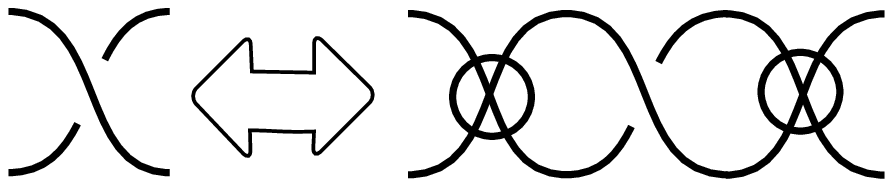}}
\caption{}\label{fig:kmove2}
\end{figure}

Kauffman \cite{rkauD} proved that the $f$-polynomial is invariant under 
the local move illustrated in Figure~\ref{fig:kmove2}, which we call {\it Kauffman's twist\/} in this paper.

\nocolon

\begin{figure}[htb]
\centerline{\epsfysize 2cm \epsfbox{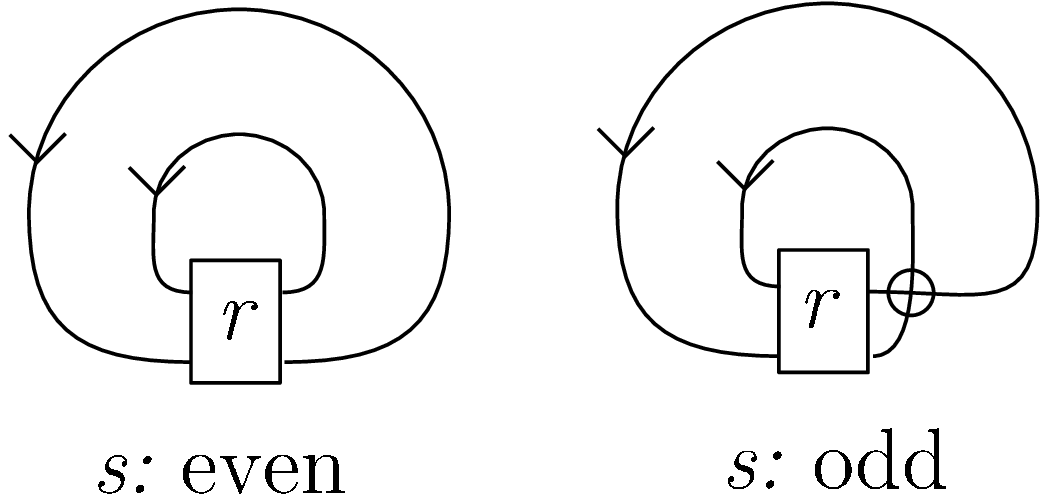}}
\caption{}\label{fig:2brknotj}
\end{figure}

Using Kauffman's twists and generalized Reidemeister moves, we see that the
$f$-polynomial of $K(r_1,\cdots, r_s)$ is equal to the $f$-polynomial of a virtual link illustrated in
Figure~\ref{fig:2brknotj}, where $r=r_1+\cdots+r_s$. If $s$ is even, then it is a $(2,r)$-torus link or a trivial
link. If
$s$ is odd and $r\ne 0$, then it is a v-alternating virtual link diagram satisfying the hypothesis of
Corollary~\ref{thm:valt3}. Thus we have the following.

\begin{cor}
\begin{itemize}
\item[\rm (1)]
If $s$ is odd and $r_1+\cdots+r_s\ne 0$, then $K(r_1,\cdots. r_s)$ is not a classical link.
\item[\rm (2)]
If $s$ is odd, $r_1+\cdots+r_s \ne 0$ and $s'$ is even, then $K(r_1,\cdots. r_s)$ and $K(r'_1,\cdots. r'_{s'})$ are
distinct virtual links. 
\end{itemize}
\end{cor}

\begin{rem}
When $s$ is even, only from a calculation of the
$f$-polynomials, we cannot conclude that $K(r_1,\cdots. r_s)$ is not a classical link. However this is true. It
will be discussed in a forthcoming paper. 
\end{rem}

\nocolon

\section{Remarks on supporting genera}\label{sec:alt}

\begin{figure}[htb]
\centerline{\epsfysize 6cm \epsfbox{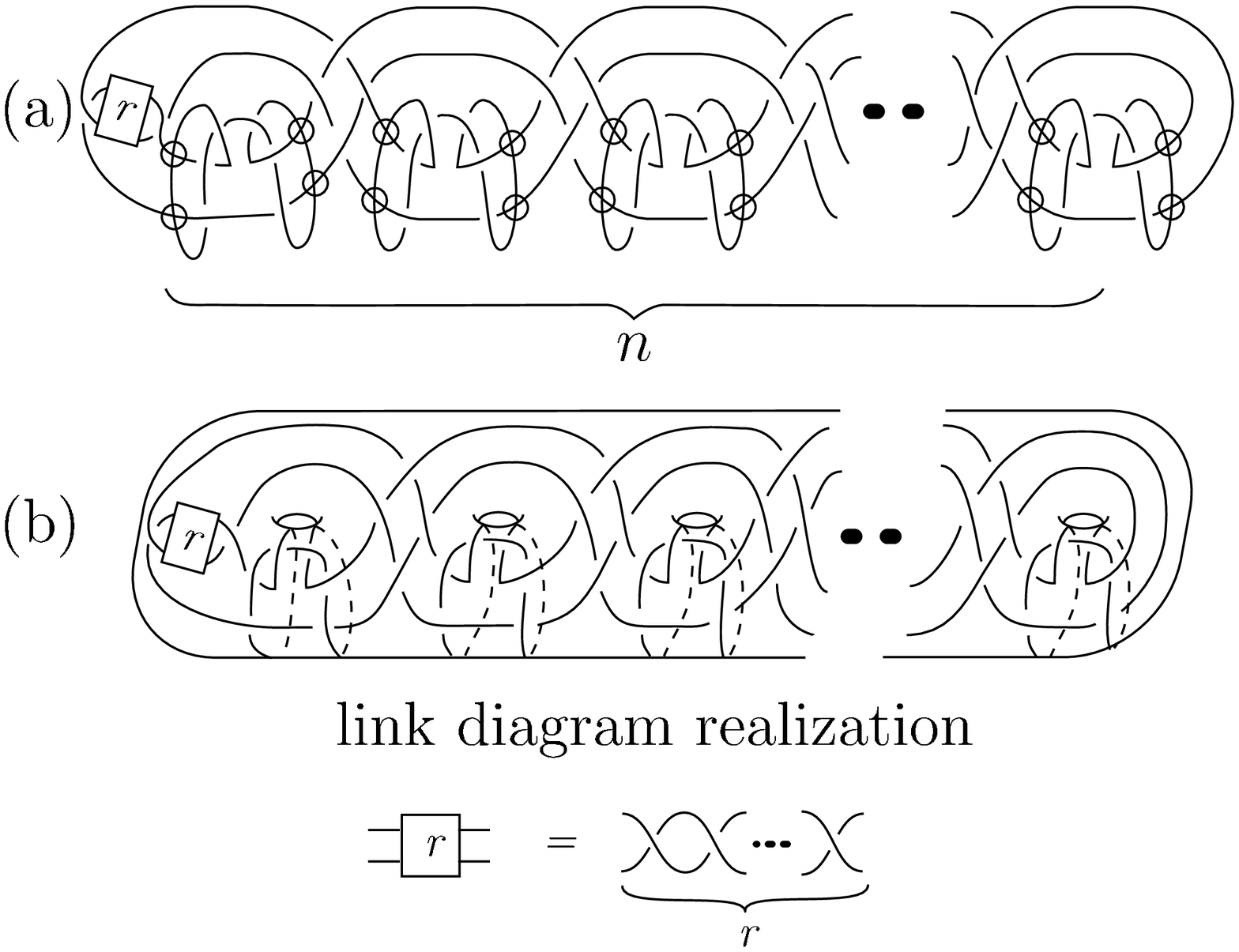}}
\caption{}\label{fig:exalt1}
\end{figure}

\begin{thm}\label{thm:altdiag1}
For any positive integer $n$, there exists an infinite family of virtual link diagrams, $D(n,r)$ $(r=0,1,2,\cdots)$,
such that
\begin{itemize}
\item[\rm (1)]
$D(n,r)$ is a proper alternating virtual link diagram, 
\item[\rm (2)]
the supporting genus is $n$, and
\item[\rm (3)]
$c(D(n,r))=10n+r-2$.
\end{itemize}
\end{thm}
\begin{proof}
A diagram $D(n,r)$ illustrated in Figure~\ref{fig:exalt1} satisfies the conditions. In the figure,  the boxed
$r$ stands for the
$r$ right half twists. The supporting genus is $n$, since it has a link diagram realization as in
Figure~\ref{fig:exalt1}(b) on a genus $n$ surface such that the complementary region consists of open disks.
\end{proof}

\begin{cor}\label{thm:altdiag2}
For any positive integer $N$, there are proper alternating (1-component) virtual link diagrams
$D_1,\cdots,D_N$ with the same crossing number and  the supporting genus of $D_k$ is $k$ ($k=1,\cdots,N$).
\end{cor}
\begin{proof}
Let $D_k$ be the diagram $D(k,10(N-k))$ introduced in Theorem~\ref{thm:altdiag1}. The crossing number of
$D_k$ is $10N-2$. 
\end{proof}

\begin{cor}\label{thm:altdiag3}
The span of the $f$-polynomial of an alternating  (1-component) virtual link $K$ is not determined only from
the number
$c(D)$ of real crossings  of a proper alternating virtual link diagram $D$ representing $K$.
\end{cor}
\begin{proof}
Let $D_1, \cdots, D_N$ be the proper alternating 1-component virtual link diagrams in the proof of
Corollary~\ref{thm:altdiag2}. Then $c(D_k)=10N-2$ and $g(D_k)=k$ for $k=1,\cdots, N$. By
Theorem~\ref{thm:alt2},
${\rm span}(D_k)=4(10N-2-k).$
Thus $D_1, \cdots, D_N$ have the same real crossing number but the spans of their $f$-polynomials are
distinct. 
\end{proof}


\nocolon

For a virtual link $L$, we define the minimal crossing number $c(L)$ and the supporting genus $g(L)$ of $L$ 
by
$$c(L)={\rm min}\{c(D)|D\text{ is a virtual link diagram representing }L\}$$
and 
$$g(L)={\rm min}\{g(D)|D\text{ is a virtual link diagram representing  }L\}.$$

In the category of classical links, the following theorem holds.

\begin{thm} [\cite{rkauA}, \cite{rmusA}, \cite{rthi}]\label{thm:kaumura2}
Let $L$ be an alternating link represented by a proper alternating link
diagram $D$. Then $c(L)=c(D)$.
\end{thm}

\begin{qst}\label{con:mincrs}
{\rm 
Let $L$ be an alternating virtual link represented by a proper alternating virtual link diagram $D$.
\begin{itemize}
\item[(1)]
Is $c(L)$  equal to $c(D)$?
\item[(2)]
Is $g(L)$ equal to $g(D)$?
\end{itemize}
By Theorem~\ref{thm:alt2}, two assertions (1) and (2)  are mutually equivalent.}
\end{qst}

As a related result, C. Adams et al. \cite{radm} and T. Kaneto \cite{rkaneB} proved the following theorem.
(C. Hayashi also informed the 
author the same result independently.)
\begin{thm}[\cite{radm}, \cite{rkaneB}]\label{thm:kane}
Let $D$ be a proper (or reduced) alternating link diagram in a closed oriented surface $F$. For any link
diagram 
$D'$ in $F$ which is related to $D$ by a finite sequence of Reidemeister moves in $F$, we have
$c(D)\leq c(D')$.
\end{thm}

This theorem is a generalization of Theorem~\ref{thm:kaumura2} when we consider that $D$ represents a link
in the thickened surface $F\times \mbox{\boldmath R}$; namely, for a link $L$ in $F\times \mbox{\boldmath
R}$ represented by a proper alternating link diagram $D$ in $F$, we have 
$c(D)=c(L),$
where $c(L)$ is the minimal crossing number of $L$ as a link in $F\times \mbox{\boldmath $R$}$.
Note that Question~\ref{con:mincrs} (1) is different from Theorem~\ref{thm:kane}.

\begin{rem}
V.O. Manturov \cite{rmanA} established another kind of generalization
of Kauffman-Murasugi-Thistlethwaite's theorem
(Theorem~\ref{thm:kaumura2}). He introduced the notion of
quasi-alternating virtual link diagram and proved that any
quasi-alternating virtual link diagram without nugatory crossing is
minimal.  A virtual link diagram is said to be {\it
quasi-alternating\/} if it is obtained from a classical alternating
link diagram by doing Kauffman's twists (Figure~\ref{fig:kmove2}) at
some crossings and virtual Reidemeister moves (in the second and third
rows of Figure~\ref{fig:rmove}).  Note that a quasi-alternating
virtual link diagram is not an alternating virtual link diagram in our
sense unless it is a classical alternating diagram or its consequences
by virtual Reidemeister moves.
\end{rem}

\rk{Acknowledgement}This research is supported by the 21st COE program
``Constitution of wide-angle mathematical basis focused on knots''.

\Addresses\recd
\end{document}